\documentclass[hidelinks]{siamart190516}

\usepackage{amsmath,amssymb}
\usepackage{smalljeff}
\pgfplotsset{compat=newest}
\usepgfplotslibrary{groupplots}

\usepgfplotslibrary{external} 
\newsiamremark{example}{Example}

\newcommand{\TheTitle}{Multivariate Rational Approximation \\
	Using a Stabilized Sanathanan-Koerner Iteration}
\newcommand{\TheAuthors}{Jeffrey M. Hokanson}
\headers{Stabilizing the SK Iteration}{\TheAuthors}

\title{{\TheTitle}\thanks{Submitted to the editors DATE.
}}

\author{Jeffrey M. Hokanson\thanks{
	Department of Computer Science, 
	University of Colorado Boulder,
	1111 Engineering Dr, Boulder, CO 80309,
	(\email{Jeffrey.Hokanson@colorado.edu}).}
}

\begin{document}
\maketitle
\begin{abstract}
	The Sanathanan-Koerner iteration developed in 1963
	is classical approach for rational approximation.
	This approach multiplies both sides of the approximation by the denominator polynomial
	yielding a linear problem and then introduces a weight at each iteration to correct for this linearization.
	Unfortunately this weight introduces a numerical instability.
	We correct this instability by constructing Vandermonde matrices
	for both the numerator and denominator polynomials
	using the Arnoldi iteration with an initial vector that enforces this weighting.
	This Stabilized Sanathanan-Koerner iteration
	corrects the instability and yields 
	accurate rational approximations of arbitrary degree.
	Using a multivariate extension of Vandermonde with Arnoldi,
	we can apply the Stabilized Sanathanan-Koerner iteration
	to multivariate rational approximation problems.
	The resulting multivariate approximations
	are often significantly better than existing techniques
	and display a more uniform accuracy throughout the domain. 
\end{abstract}

\begin{keywords}
multivariate rational approximation,
Sanathanan-Koerner iteration,
Vandermonde with Arnoldi,
least squares
\end{keywords}
\begin{AMS}
%	30E10, % Approximation in the complex domain
	41A20, % Approximation by rational functions (used in AAA paper)
	41A63, % Multidimensional problems, used in Anthony's paper AKL+19x
%	26C15, % Rational functions
	65D15 % Algorithms for approximation of functions (used in AAA paper)
%	30E05, % Moment problems, interpolation problems
%	62J02, %  	General nonlinear regression
%	90C53 % mathematical programming,   	Methods of quasi-Newton type
\end{AMS}
\begin{DOI}

\end{DOI}

% Two issues
% -> stable basis:
%   - cite Berruit 1988 for barycentric rational form 
% -> spurious/unfortunate/large-residual local minimizers

\section{Introduction\label{sec:intro}}
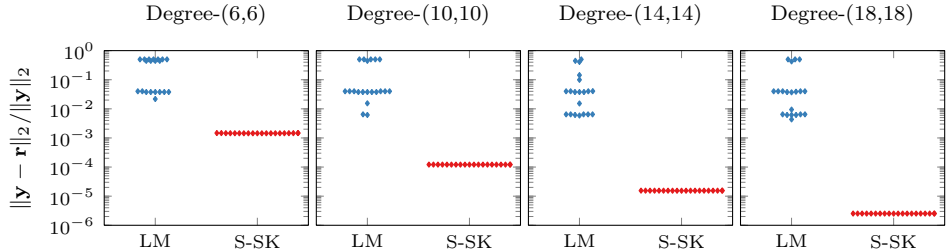
\begin{figure}

\begin{tikzpicture}
\tikzset{mark size = 1pt}
\begin{groupplot}[
	group style = {group size = 4 by 1, horizontal sep = 3pt},
	width = 0.33\linewidth,
	height = 0.3\linewidth,
	ymode = log,
	ymin = 1e-6, 
	ymax = 1e0,
	ytickten = {-6, -5, ..., 1},
	xmin = -0.5, xmax = 1.5,
	xtick = {0,1},
	xticklabels = {LM, S-SK},
]
	\nextgroupplot[
		title = {Degree-(6,6)},
		ylabel = {$\|\ve y- \ve r\|_2/\|\ve y\|_2$},
	]

	\addplot[blue, only marks, mark=diamond*]
		table [x=x, y expr= (10^\thisrow{y})/258.20018074483727] {data/fig_init_06_opt.dat}; 
%	\addplot[green, only marks, mark=*, mark size = 1.5pt]
%		table [x=x, y expr= 10^\thisrow{y}] {data/fig_init_10_vf.dat}; 
	\addplot[red, only marks, mark=diamond*]
		table [x expr=\thisrow{x}-1, y expr= (10^\thisrow{y})/258.20018074483727] {data/fig_init_06_ssk.dat}; 
	
	\nextgroupplot[
		title = {Degree-(10,10)},
		yticklabels = {,,},
	]

	\addplot[blue, only marks, mark=diamond*]
		table [x=x, y expr= (10^\thisrow{y})/258.20018074483727] {data/fig_init_10_opt.dat}; 
%	\addplot[green, only marks, mark=*, mark size = 1.5pt]
%		table [x=x, y expr= 10^\thisrow{y}] {data/fig_init_10_vf.dat}; 
	\addplot[red, only marks, mark=diamond*]
		table [x expr=\thisrow{x}-1, y expr= (10^\thisrow{y})/258.20018074483727] {data/fig_init_10_ssk.dat}; 
	
	\nextgroupplot[
		title = {Degree-(14,14)},
		yticklabels = {,,},
	]
	\addplot[blue, only marks, mark=diamond*]
		table [x =x, y expr= (10^\thisrow{y})/258.20018074483727] {data/fig_init_14_opt.dat}; 
%	\addplot[green, only marks, mark=*]
%		table [x=x, y expr= 10^\thisrow{y}] {data/fig_init_14_vf.dat}; 
	\addplot[red, only marks, mark=diamond*]
		table [x expr=\thisrow{x}-1, y expr= (10^\thisrow{y})/258.20018074483727] {data/fig_init_14_ssk.dat}; 
	
	\nextgroupplot[
		title = {Degree-(18,18)},
		yticklabels = {,,},
	]
	\addplot[blue, only marks, mark=diamond*]
		table [x=x, y expr= (10^\thisrow{y})/258.20018074483727] {data/fig_init_18_opt.dat}; 
%	\addplot[green, only marks, mark=*]
%		table [x=x, y expr= 10^\thisrow{y}] {data/fig_init_18_vf.dat}; 
	\addplot[red, only marks, mark=diamond*]
		table [x expr=\thisrow{x}-1, y expr= (10^\thisrow{y})/258.20018074483727] {data/fig_init_18_ssk.dat}; 

\end{groupplot}
\end{tikzpicture}
\caption{The residual norm of the rational approximation
from twenty different initializations
computed using Levenberg-Marquardt (LM) applied to~\cref{eq:ls_rat}
and the Stabilized SK iteration (S-SK) introduced here.
In this example, we approximate $|x|$ on $[-1,1]$ using 200,000 equispaced spaced points.
The initial numerator and denominator polynomials, both of degree $m$,
take normally distributed values with zero mean and unit variance
on $m+1$ equispaced points on $[-1,1]$.
}
\label{fig:init}
\end{figure}

Given pairs of inputs $\lbrace x_j \rbrace_{j=1}^M \subset \C$
and outputs $\lbrace y_j \rbrace_{j=1}^M\subset \C$,
we wish to construct a degree-$(m,n)$ rational approximation $r:\C\to\C$
where
\begin{equation}
	y_j \approx r_j  = r( x_j) \coloneqq \frac{ p( x_j)}{q( x_j)}
\end{equation}
and $p$ and $q$ are two polynomials of degree $m$ and $n$,
denoted $p \in \set P_m$ and $q \in \set P_n$.
After constructing discrete bases $\ma P\in \C^{M\times(m+1)}$ and $\ma Q\in \C^{M\times(n+1)}$ 
for $\set P_m$ and $\set P_n$ on $\lbrace x_j \rbrace_{j=1}^M$,
we can restate the rational approximation problem 
as identifying polynomial coefficients $\ve a\in \C^{m+1}$ and $\ve b \in \C^{n+1}$
such that 
\begin{equation}\label{eq:rat_approx}
	\ve y \approx \diag(\ma Q \ve b)^{-1} \ma P \ve a.
\end{equation}
One challenge of rational approximation is
that as a nonlinear least squares problem,
%Posed in this form,
%we must construct well-conditioned discrete polynomial bases on $\lbrace x_j \rbrace_{j=1}^M$
%for $\set P_m$ and $\set P_n$.
\begin{equation}\label{eq:ls_rat}
	\min_{\ve a \in \C^{m+1}, \ve b \in \C^{n+1}}
		\| \ve y - \diag(\ma Q \ve b)^{-1} \ma P \ve a\|_2,
\end{equation}
most solvers when initialized randomly tend to converge
to approximations with a large residual norm as illustrated in \cref{fig:init}.
This had lead to a variety of non-optimal techniques based on 
\emph{linearizing} the rational approximation problem
by multiplying both sides of~\cref{eq:rat_approx} by the denominator:
\begin{equation}\label{eq:linearized}
	\diag(\ma Q \ve b)^{-1} \ma P \ve a \approx \ve y
	\quad \Longrightarrow \quad
	\ma P \ve a \approx \diag(\ma Q \ve b) \ve y =  \diag(\ve y) \ma Q \ve b.
\end{equation}
Rational approximation algorithms using this linearization, 
although not least-squares optimal in general~\cite{Whi87},
tend to yield rational approximations with smaller residual norm.
These include:
\emph{linearized rational approximation}
that solves~\cref{eq:linearized} in a least-squares sense~\cite{AKL+19x, Lev59},
the \emph{Sanathanan-Koerner} (SK) \emph{iteration}~\cite{SK63},
\emph{Vector Fitting}~\cite{Gus06,GS99},
the \emph{Loewner framework}~\cite{AA86},
and \emph{Adaptive Anderson-Antoulas} (AAA)~\cite{NST18}.
An important consideration in these algorithms 
is the choice polynomial basis to construct $\ma P$ and $\ma Q$.
For example, Vector Fitting uses a barycentric Lagrange basis~\cite{BT04}
and iteratively updates the interpolation nodes
whereas AAA uses the same basis but adds nodes greedily.
Here we construct a well-conditioned basis for the SK iteration
using a weighted Arnoldi iteration.

\subsection{The SK Iteration}
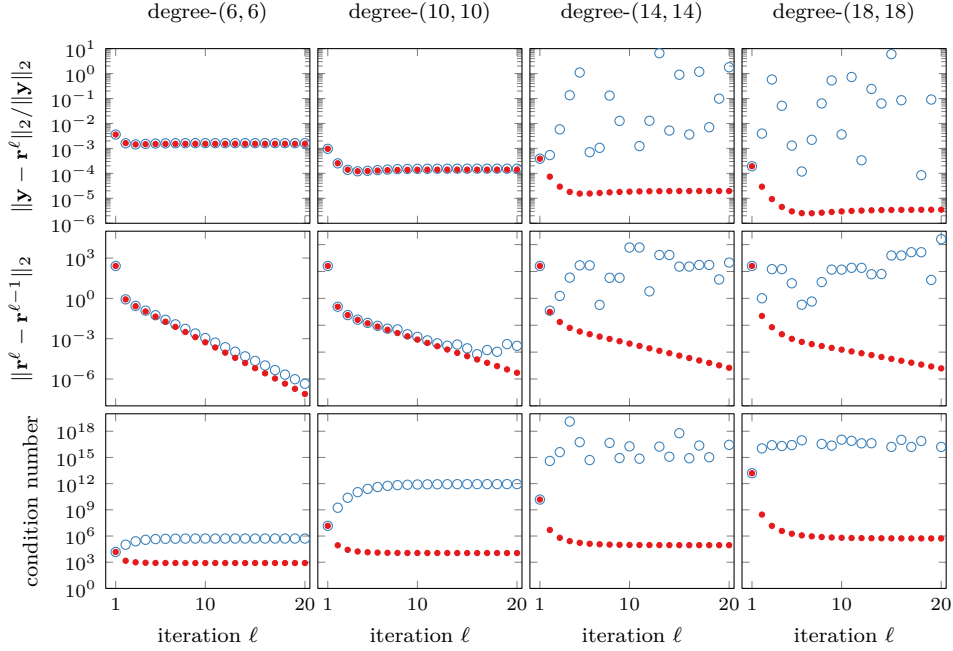
\begin{figure}
\centering
\begin{tikzpicture}
\begin{groupplot}[
	group style = {group size = 4 by 3, vertical sep = 3pt, horizontal sep = 3pt},
	ymode = log,
	xmin = 0,
	xmax = 20.5,
	ymin = 1e-12,
	ymax = 1e1,
	width = 0.33\linewidth,
	height = 0.3\linewidth,
	xtick = {1,10,20},
	ytickten = {-12,-10,...,20},
	every axis plot post/.append style={
		mark size = 1pt,
		},
	]
	\nextgroupplot[
		title = {degree-$(6,6)$},
		ylabel = {$\|\ve y - \ve r^\ell\|_2/\|\ve y\|_2$},
		xticklabels = {,,},
		ymin = 1e-6, ymax = 1e1,
		ytickten = {-6,-5,...,1},
		]
		\addplot[blue, only marks, mark=o, mark options ={scale=1.75}]
			table [x expr=\coordindex+1, y expr = \thisrow{err}/258.2] {data/fig_abs_arnoldi_6_6.dat};
		\addplot[red, only marks, mark=*]
			table [x expr=\coordindex+1, y expr = \thisrow{err}/258.2] {data/fig_abs_arnoldi_rebase_6_6.dat};
	
	\nextgroupplot[
		title = {degree-$(10,10)$},
		yticklabels = {,,},
		xticklabels = {,,},
		ymin = 1e-6, ymax = 1e1,
		ytickten = {-6,-5,...,1},
		]
		\addplot[blue, only marks, mark=o, mark options = {scale =1.75}]
			table [x expr=\coordindex+1, y expr = \thisrow{err}/258.2] {data/fig_abs_arnoldi_10_10.dat};
		\addplot[red, only marks, mark=*]
			table [x expr=\coordindex+1, y expr = \thisrow{err}/258.2] {data/fig_abs_arnoldi_rebase_10_10.dat};
		
	\nextgroupplot[
		title = {degree-$(14,14)$},
		yticklabels = {,,},
		xticklabels = {,,},
		ymin = 1e-6, ymax = 1e1,
		ytickten = {-6,-5,...,1},
		]
		\addplot[blue, only marks, mark=o, mark options = {scale = 1.75}]
			table [x expr=\coordindex+1, y expr= \thisrow{err}/258.2] {data/fig_abs_arnoldi_14_14.dat};
		\addplot[red, only marks, mark=*]
			table [x expr=\coordindex+1, y expr = \thisrow{err}/258.2] {data/fig_abs_arnoldi_rebase_14_14.dat};

	\nextgroupplot[
		title = {degree-$(18,18)$},
		yticklabels = {,,},
		xticklabels = {,,},
		ymin = 1e-6, ymax = 1e1,
		ytickten = {-6,-5,...,1},
		]
		\addplot[blue, only marks, mark=o, mark options = {scale = 1.75}]
			table [x expr=\coordindex+1, y expr= \thisrow{err}/258.2] {data/fig_abs_arnoldi_18_18.dat};
		\addplot[red, only marks, mark=*]
			table [x expr=\coordindex+1, y expr= \thisrow{err}/258.2] {data/fig_abs_arnoldi_rebase_18_18.dat};

	%%%%%%%% Second row %%%%%%%%%%%%%%%%%%%%%%%%%%%%%%%%%%%%%
	\nextgroupplot[
		ylabel = {$\|\ve r^\ell - \ve r^{\ell-1}\|_2$},
		xticklabels = {,,},
		ymin = 1e-8, ymax = 1e5,
		ytickten = {-9,-6,...,6},
		]
		\addplot[blue, only marks, mark=o, mark options ={scale=1.75}]
			table [x expr=\coordindex+1, y = delta_fit] {data/fig_abs_arnoldi_6_6.dat};
		\addplot[red, only marks, mark=*]
			table [x expr=\coordindex+1, y = delta_fit] {data/fig_abs_arnoldi_rebase_6_6.dat};
	
	\nextgroupplot[
		yticklabels = {,,},
		xticklabels = {,,},
		ymin = 1e-8, ymax = 1e5,
		]
		\addplot[blue, only marks, mark=o, mark options = {scale =1.75}]
			table [x expr=\coordindex+1, y = delta_fit] {data/fig_abs_arnoldi_10_10.dat};
		\addplot[red, only marks, mark=*]
			table [x expr=\coordindex+1, y = delta_fit] {data/fig_abs_arnoldi_rebase_10_10.dat};
		
	\nextgroupplot[
		yticklabels = {,,},
		xticklabels = {,,},
		ymin = 1e-8, ymax = 1e5,
		]
		\addplot[blue, only marks, mark=o, mark options = {scale = 1.75}]
			table [x expr=\coordindex+1, y = delta_fit] {data/fig_abs_arnoldi_14_14.dat};
		\addplot[red, only marks, mark=*]
			table [x expr=\coordindex+1, y = delta_fit] {data/fig_abs_arnoldi_rebase_14_14.dat};

	\nextgroupplot[
		yticklabels = {,,},
		xticklabels = {,,},
		ymin = 1e-8, ymax = 1e5,
		]
		\addplot[blue, only marks, mark=o, mark options = {scale = 1.75}]
			table [x expr=\coordindex+1, y = delta_fit] {data/fig_abs_arnoldi_18_18.dat};
		\addplot[red, only marks, mark=*]
			table [x expr=\coordindex+1, y = delta_fit] {data/fig_abs_arnoldi_rebase_18_18.dat};

	%%%%%%%% Third row %%%%%%%%%%%%%%%%%%%%%%%%%%%%%%%%%%%%%
	\nextgroupplot[
		xlabel = {iteration $\ell$},
		ylabel = {condition number},
		ymin = 1e0, ymax = 1e20,
		ytickten = {0,3,...,21},
		]
		\addplot[blue, only marks, mark=o, mark options ={scale=1.75}]
			table [x expr=\coordindex+1, y = cond] {data/fig_abs_arnoldi_6_6.dat};
		\addplot[red, only marks, mark=*]
			table [x expr=\coordindex+1, y = cond] {data/fig_abs_arnoldi_rebase_6_6.dat};
	
	\nextgroupplot[
		xlabel = {iteration $\ell$},
		yticklabels = {,,},
		ymin = 1e0, ymax = 1e20,
		ytickten = {0,3,...,21},
		]
		\addplot[blue, only marks, mark=o, mark options = {scale =1.75}]
			table [x expr=\coordindex+1, y = cond] {data/fig_abs_arnoldi_10_10.dat};
		\addplot[red, only marks, mark=*]
			table [x expr=\coordindex+1, y = cond] {data/fig_abs_arnoldi_rebase_10_10.dat};
		
	\nextgroupplot[
		xlabel = {iteration $\ell$},
		yticklabels = {,,},
		ymin = 1e0, ymax = 1e20,
		ytickten = {0,3,...,21},
		]
		\addplot[blue, only marks, mark=o, mark options = {scale = 1.75}]
			table [x expr=\coordindex+1, y = cond] {data/fig_abs_arnoldi_14_14.dat};
		\addplot[red, only marks, mark=*]
			table [x expr=\coordindex+1, y = cond] {data/fig_abs_arnoldi_rebase_14_14.dat};

	\nextgroupplot[
		xlabel = {iteration $\ell$},
		yticklabels = {,,},
		ymin = 1e0, ymax = 1e20,
		ytickten = {0,3,...,21},
		]
		\addplot[blue, only marks, mark=o, mark options = {scale = 1.75}]
			table [x expr=\coordindex+1, y = cond] {data/fig_abs_arnoldi_18_18.dat};
		\addplot[red, only marks, mark=*]
			table [x expr=\coordindex+1, y = cond] {data/fig_abs_arnoldi_rebase_18_18.dat};

\end{groupplot}
\end{tikzpicture}
\caption{Iteration histories for approximating $|x|$ on $[-1,1]$ using 200,000 equispaced points
using the SK iteration.
Hollow dots show the standard SK iteration using the Arnoldi basis with condition number one;
the solid dots show the Stabilized SK iteration;
$\ve r^{\ell}$ is the rational approximation evaluated at $\lbrace x_j \rbrace_{j=1}^M$ on the $\ell$th iteration;
the bottom row refers to the condition number of the smallest right singular vector
of the system in~\cref{eq:skls} or~\cref{eq:sskls}; see, e.g.,~\cite[Chap.~3, eq.~(3.16)]{Ste01}.
}
\label{fig:abs}
\end{figure}

Sanathanan and Koerner's key contribution was to introduce
a weight into the linearized rational approximation problem~\cref{eq:linearized}
to better reflect the original rational approximation problem~\cref{eq:rat_approx}.
At the $\ell+1$th iteration, 
they include the weight $\diag(\ma Q \ve b^{\ell})^{-1}$
(the previous iterate's denominator)
and compute new coefficients $\ve a^{\ell+1}$ and $\ve b^{\ell+1}$ 
solving the approximation problem
\begin{equation}~\label{eq:skapprox}
	\diag(\ma Q \ve b^{\ell})^{-1} \ma P \ve a^{\ell+1} \approx 
	\diag(\ma Q \ve b^{\ell})^{-1} \diag(\ve y) \ma Q \ve b^{\ell+1}.
\end{equation}
If $\ve a^\ell \to \ve a^\star$ and $\ve b^\ell \to \ve b^\star$,
then this limit is a rational approximation of $\ve y$:
\begin{equation}
	\diag(\ma Q \ve b^\star)^{-1} \ma P \ma a^\star \approx \ve y. 
\end{equation}
If we approximate in a least-squares sense,
solving~\cref{eq:skapprox} corresponds to 
\begin{equation}\label{eq:skls}
	\ve a^{\ell+1}, \ve b^{\ell+1} \leftarrow
		\argmin_{\substack{\ve a, \ve b \\ \ve b \ne \ve 0}}
			\left\| \diag(\ma Q \ve b^{\ell} )^{-1} 
				\begin{bmatrix}  \ma P  & - \diag(\ve y) \ma Q \end{bmatrix}
				\begin{bmatrix} \ve a \\ \ve b \end{bmatrix} 
			\right\|_2
\end{equation}
which we solve using the singular value decomposition (SVD).
The primary challenge with the SK iteration is even if $\ma P$ and $\ma Q$ have condition number one,
the weighting can cause this system to be ill-conditioned
and consequently yield poor approximations as illustrated in \cref{fig:abs}.

\subsection{Stabilizing the SK Iteration}
We correct the ill-conditioning of the SK iteration
by building $\ma P$ and $\ma Q$ using a weighted Arnoldi iteration.
Recall the Arnoldi iteration builds an orthonormal basis for the 
for the Krylov subspace
\begin{equation}
	\set K_{m}(\ma A, \ve w) = 
		\Span \left \lbrace
			\ve w, \ma A \ve w, \ma A^2 \ve w, \ldots, \ma A^{m-1} \ve w 
		\right \rbrace
\end{equation}
by applying Gram-Schmidt to produce orthonormal vectors $\ve q_\ell$
from the sequence $\ve v_1 = \ve w$ and $\ve v_\ell \leftarrow \ma A  \ve q_{\ell-1}$.
If use the Arnoldi iteration to construct a basis for 
\begin{equation}
	\set K_{m+1}(\diag (\ve x), \ve 1) =
			\Span \left \lbrace 
				\begin{bmatrix} 1 \\ \vdots \\ 1  \end{bmatrix},
				\begin{bmatrix} x_1 \\ \vdots \\ x_M  \end{bmatrix},
				\begin{bmatrix} x_1^2 \\ \vdots \\ x_M^2  \end{bmatrix},
				\cdots,
				\begin{bmatrix} x_1^m \\ \vdots \\ x_M^m  \end{bmatrix}
			\right\rbrace 
\end{equation}
we have constructed an orthonormal basis on $\lbrace x_j \rbrace_{j=1}^M$ 
for polynomials of degree $m$.
This technique is called \emph{Vandermonde with Arnoldi}~\cite{BNT19x}
and accurately computes a discrete polynomial basis while avoiding
the ill-conditioning of standard Vandermonde matrices~\cite{Pan16}.
% mention numerical stability due to multiplies?
Here we use this insight to construct an orthonormal basis with respect to
weighting at each step of the SK iteration.
If $\ma V_m \in \C^{M\times (m+1)}$ is a basis for $\set P_m$ on $\lbrace x_j \rbrace_{j=1}^M$, 
then $\Range(\ma V_m) = \set K_{m+1}(\diag(\ve x), \ve 1)$ and 
\begin{equation}
	\Range (\diag(\ve w) \ma V_m)
	= \Span \left \lbrace
				\diag(\ve w) \begin{bmatrix} x_1^k \\ \vdots \\  x_M^k  \end{bmatrix}
			\right \rbrace_{k=0}^m
	= \set K_{m+1}(\diag(\ve x), \ve w).
\end{equation}
By choosing the initial vector $\ve w^\ell$
to be the inverse of the previous iterate's denominator, $w_j^\ell = 1/q_{\ell-1}(x_j)$,
we can then construct iteration-dependent bases
$\ma P^\ell$ and $\ma Q^\ell$ using Vandermonde with Arnoldi.
As these bases implicitly include the weight,
we can update the polynomial coefficients by solving
\begin{equation}\label{eq:sskls}
	\ve a^\ell, \ve b^\ell \leftarrow \min_{\substack{\ve a, \ve b\\ \ve b\ne 0}} 
		\left\| 
			\begin{bmatrix} \ma P^\ell & -\diag(\ve y) \ma Q^\ell \end{bmatrix}
			\begin{bmatrix} \ve a \\ \ve b \end{bmatrix}
		\right\|_2.
\end{equation}
Unlike the standard SK iteration~\cref{eq:skls},
this problem tends to be well-conditioned  
and often converges linearly to its fixed points as seen in \cref{fig:abs}.
As the same weight appears in both the numerator and denominator,
we can evaluate the rational approximation on $\lbrace x_j\rbrace_{j=1}^M$ by simply computing 
$\diag(\ma Q^\ell \ve b^\ell)^{-1}\ma P^\ell \ve a^\ell$.
However to evaluate the denominator when computing the initial vector $\ve w^\ell$,
we must undo the action of the previous weight:
$q_{\ell-1}(x_j) = [\ma Q^{\ell-1}\ve b^{\ell-1}]_j/w_j^{\ell-1}$.
This new \emph{Stabilized Sanathanan-Koerner iteration} is summarized in \cref{alg:isk}.

\begin{algorithm}[t]
\begin{minipage}{\linewidth}
\begin{algorithm2e}[H]
\Input{Data $\lbrace x_j, y_j \rbrace_{j=1}^M$, degrees $m$, $n$}
\Output{$\ma R_\ma P^\ell$, $\ma R_\ma Q^\ell$, $\ve a^\ell$, $\ve b^\ell$
	for $\ell$
	minimizing $\| \ve y - \diag(\ma Q^\ell \ve b^\ell)^{-1}\ma P^\ell \ve a^\ell\|_2$}
$\ve w^0 \leftarrow \ve 1$\;
\For{$\ell=0,1,2,\ldots$ and not converged}{
	$\ma P^\ell, \ma R_\ma P^\ell \leftarrow$ Arnoldi for $\set K_{m+1}(\diag(\ve x), \ve w^\ell)$
		or \cref{alg:vandermonde} if multivariate%
		\;\label{alg:isk:P}
	$\ma Q^\ell, \ma R_\ma Q^\ell \leftarrow$ Arnoldi for $\set K_{n+1}(\diag(\ve x), \ve w^\ell)$
		or \cref{alg:vandermonde} if multivariate%
		\;\label{alg:isk:Q}
	%$\ve v \leftarrow$ smallest right singular value of 
	%	$\begin{bmatrix} \ma P^\ell & -\diag(\ve y) \ma Q^{\ell}\end{bmatrix}$\;
	%$\ve a^\ell \leftarrow \ve v[0:m+1], \quad \ve b^\ell \leftarrow \ve v[m+2:m+n+2]$\;
	$\ve a^\ell, \ve b^\ell \leftarrow 
		\min_{\ve a, \ve b} \| \ma P^\ell \ve a - \diag(\ve y) \ma Q^\ell \ve b \|_2$
		s.t. $\|\ve a\|_2^2 + \|\ve b \|_2^2 = 1$\;
	%$\ve w^{\ell+1} \leftarrow \diag[\diag(\ve w^\ell)^{-1}\ma Q^\ell \ve b^{\ell}]^{-1} \ve 1 $\;
	$[\ve w^{\ell+1}]_j \leftarrow w_j^{\ell}/[\ma Q^\ell \ve b^\ell]_j$\;
	\label{alg:isk:weight}
}

\end{algorithm2e}
\vspace{-1.5em}
\end{minipage}
\caption{Stabilized Sanathanan-Koerner Iteration}
\label{alg:isk}
\end{algorithm}

\subsection{Advantages}
The main utility of the Stabilized SK iteration
comes in its use for multivariate rational approximation.
The only modification required is to replace the use of Vandermonde with Arnoldi
to construct discrete polynomial bases $\ma P^\ell$ and $\ma Q^\ell$
with the corresponding multivariate generalization developed in~\cite[Subsec.~3.2]{AKL+19x}.
The rational approximations generated by the Stabilized SK iteration 
often have a least-squares residual norm an order of magnitude smaller
than both \emph{Parametric-AAA} (p-AAA)~\cite{CG20x} 
and the linearized approach advocated in~\cite[Subsec.~3.1]{AKL+19x};
moreover the Stabilized SK iteration avoids the spurious poles 
often encountered in other algorithms (\cref{sec:examples:spurious}).
Additionally, the Stabilized SK places no restriction on the points $\lbrace \ve x_j \rbrace_{j=1}^M\subset \C^d$,
unlike p-AAA which requires points on a tensor-product grid.

Applied to univariate rational approximation problems,
the Stabilized SK iteration yields comparable approximations to 
Vector Fitting (\cref{sec:examples:univariate}).
These approximations often have a far smaller least-squares residual norm
than those generated by AAA and the linearized approach.

\subsection{Disadvantages}
Unfortunately the Stabilized SK iteration inherits some limitations
of the original SK iteration:
fixed points of this iteration are not least squares optimal~\cite[Subsec.~5.2]{Whi87},
the iteration can cycle (this happens for odd degrees in the example from \cref{fig:abs}),
and iterates do not monotonically decrease the least-squares residual norm as seen in \cref{fig:abs}.
We address the last two issues by performing only a few iterations (typically twenty)
and returning the best rational approximation.
The first issue tends not to be significant in practice.
As residual of the rational approximation decreases,
fixed points of the SK iteration approach those of the least squares problem.
Often we see in our examples that refining rational approximation using nonlinear least squares
only slightly decreases the residual norm.

The Stabilized SK iteration is more expensive than other algorithms
due to the need to perform two orthogonalizations at each step.
However, this does not increase the asymptotic complexity;
each of AAA, SK, Stabilized SK, and Vector Fitting 
require $\order(MN^2)$ operations
where $N$ is the number of columns in $\ma P$ and $\ma Q$.

\subsection{Outline}
In the remainder of this paper we
first review the multivariate Vandermonde with Arnoldi algorithm introduced in~\cite{AKL+19x}
and extend it for total-degree polynomials (\cref{sec:vandermonde}).
We then briefly discuss implementation details 
for refining rational approximations using nonlinear least squares techniques (\cref{sec:refinement}).
Finally, we conclude with several numerical examples comparing 
the Stabilized SK iteration to other rational approximation techniques
on both univariate and multivariate test problems (\cref{sec:examples}).

\subsection{Reproducibility}
Following the principles of reproducible research,
we provide software implementing the algorithms in
this paper and scripts generating the figures
at {\tt \url{https://github.com/jeffrey-hokanson/polyrat}}.

\section{Multivariate Vandermonde with Arnoldi\label{sec:vandermonde}}
We extend the univariate Stabilized Sanathanan-Koerner iteration
to multivariate rational approximation,
\begin{equation}
	y_j \approx r_j = r(\ve x_j) \coloneqq \frac{p(\ve x_j)}{q(\ve x_j)}
		\quad \text{where} \quad
		\lbrace \ve x_j \rbrace_{j=1}^M \subset \C^d
		\quad \text{and} \quad
		p,q \text{ polynomials},
\end{equation}
by replacing Vandermonde with Arnoldi 
%on lines~\ref{alg:isk:P}-\ref{alg:isk:Q} of \cref{alg:isk}
with its multivariate extension developed in~\cite[Subsec.~3.2]{AKL+19x}.
Here consider two classes of multivariate polynomials:
\emph{total degree polynomials} $\set P^{\text{tot}}_m$
and \emph{maximum degree polynomials} $\set P_{\ve m}^{\text{max}}$,
\begin{align}
	\label{eq:tot}
	\set P_{m}^{\text{tot}}
		& \coloneqq \Span \left\lbrace f:
			f(\ve x) = \prod_{i=1}^d x_i^{\alpha_i}
			\right\rbrace_{|\ve \alpha| \le m},
		& \text{where} & \quad |\ve \alpha| = \sum_{i=1}^d \alpha_i; \\
	\label{eq:max}
	\set P_{\ve m}^{\text{max}}
		& \coloneqq 
			\Span \left \lbrace f:
				f(\ve x) = \prod_{i=1}^d x_i^{\alpha_i}
			\right\rbrace_{ \ve \alpha \le \ve m},
		& \text{where} & \quad \ve \alpha \le \ve m \Leftrightarrow
			\alpha_i \le m_i.
\end{align}
The main difference in the multivariate extension of Vandermonde with Arnoldi
is that the basis no longer corresponds to a Krylov subspace.
Instead we generate new columns by carefully selecting
one coordinate of the points $\lbrace \ve x_j \rbrace_{j=1}^M$ to multiply by a proceeding column
and then apply Gram-Schmidt as before.
Here we briefly provide the details on constructing this basis
and evaluating this basis at new points.

\subsection{Building a Basis}
Our goal will be to find an ordering of multi-indices $\set I$
appearing in the polynomial basis definition in~\cref{eq:tot,eq:max}
such that the columns
$\ve q_1, \ldots, \ve q_\ell$ generated by multivariate Vandermonde with Arnoldi satisfy 
\begin{equation}\label{eq:mv_constraint}
	\Span \lbrace \ve q_k \rbrace_{k=1}^\ell
	= \Span \left \lbrace
		\begin{bsmallmatrix}
			\ve x_1^{\ve \alpha} \\ \vdots \\ \ve x_M^{\ve \alpha}
		\end{bsmallmatrix}
	\right\rbrace_{\ve \alpha \in \set I[1:\ell]},
	\quad \text{where} \quad 
	\ve x^{\ve \alpha} = \prod_{i=1}^d x_i^{\alpha_i}.
\end{equation}
At each step of multivariate Vandermonde with Arnoldi,
we generate the next column $\ve v_{\ell+1}$ 
to orthogonalize against $\ve q_1, \ldots, \ve q_\ell$ 
by multiplying $\ve q_k$ by the $j$th:
\begin{equation}\label{eq:mva_update}
	\ve v_{\ell+1} = \diag\left( 
			\begin{bmatrix} [\ve x_1]_j \\ \vdots \\ [\ve x_M]_j \end{bmatrix}
			\right)
			\ve q_k. 
\end{equation}
We choose $j$ and $k$ by finding the smallest $k$ such that 
$\set I[k] + \ve e_j = \set I[\ell+1]$
where $\ve e_j$ is the $j$th column of the $d\times d$ identity matrix.
With this update rule, we need to pick an ordering $\set I$ 
such that~\cref{eq:mv_constraint} is satisfied;
many orderings do not satisfy this constraint!
For total degree polynomials a \emph{grevlex} ordering
(ordered by total degree and then lexicographically)
satisfies this constraint; i.e.,
\begin{equation*}
	\set I_{3}^{\text{tot}} = \left[
		(0,0), (1, 0), (0,1), (2,0), (1,1), (0,2)
		,(3,0), (2,1), (1,2), (0,3)
	\right].
\end{equation*}
For maximum degree polynomials we satisfy this constraint
using a lexicographic ordering; i.e., 
\begin{equation*}
	\set I_{(2,2)}^{\text{max}} = \left[
		(0,0), (0,1), (0,2),
		(1,0), (1,1), (1,2),
		(2,0), (2,1), (2,2)
	\right].
\end{equation*}
\Cref{alg:vandermonde} summarizes the multivariate Vandermonde with Arnoldi process.
In our implementation we use classical Gram-Schmidt with two steps of iterative refinement~\cite[Sec.~6]{Bjo94}
rather than modified Gram-Schmidt.
Although this uses more floating point operations,
classical Gram-Schmidt allows us to make use of BLAS level 2 operations
yielding a net decrease in wall-clock time compared to 
the BLAS level 1 operations used in modified Gram-Schmidt.

\begin{algorithm}[t]
\begin{minipage}{\linewidth}
\begin{algorithm2e}[H]
	\Input{ $\lbrace \ve x_j \rbrace_{j=1}^M \subset \C^d$, weight $\ve w$ 
			index set $\set I$ of length $N$}
	\Output{ $\ma Q\in \C^{M\times N}$, $\ma R \in \C^{N\times N}$}
	$\ma Q \leftarrow \ma 0$, $\ma R \leftarrow \ma 0$\;
	$[\ma R]_{1,1} \leftarrow \|\ve w\|_2$\;
	$[\ma Q]_{:,1} \leftarrow \ve w /[\ma R]_{1,1}$\;
	\For{$\ell=2,3, \ldots, |\set I|$}{
		Pick smallest $k$ such that $\exists j$ where  $\set I[k] + \ve e_j = \set I[\ell]$\;
		$\ve v_\ell \leftarrow \diag([[\ve x_1]_j, \cdots, [\ve x_M]_j]) \ve q_k$\;
		\For{$t=1,2$}{
			$\ve s \leftarrow [\ma Q]_{\cdot,1:\ell-1}^* \ve v_\ell$\;
			$\ve v_\ell \leftarrow \ve v_\ell - [\ma Q]_{\cdot,1:\ell-1} \ve s$\;
			$[\ma R]_{1:\ell-1,\ell} \leftarrow [\ma R]_{1:k-1,k} + \ve s$\;
		}
		$[\ma R]_{\ell,\ell} \leftarrow \|\ve v_\ell\|_2$\;
		$[\ma Q]_{:,\ell} \leftarrow \ve v_\ell/[\ma R]_{\ell,\ell}$\;
	}
	\vspace*{-1.5em}
\end{algorithm2e}
\end{minipage}
\caption{Multivariate Vandermonde with Arnoldi}
\label{alg:vandermonde}
\end{algorithm}

\subsection{Evaluating a Basis}
Once we have constructed a polynomial basis in \cref{alg:vandermonde},
we need to be able to evaluate the resulting basis at new points $\ve z\in \C^d$.
To do so, we simply repeat the construction as before but keep $\ma R$ fixed
as illustrated in \cref{alg:arnoldi_eval}.

\begin{algorithm}[t]
\begin{minipage}{\linewidth}
\begin{algorithm2e}[H]
	\Input{ $\lbrace \ve z_j \rbrace_{j=1}^M \subset \C^d$, 
			index set $\set I$ of length $N$, $\ma R\in \C^{N\times N}$ }
	\Output{ $\ma W\in \C^{M\times N}$}
	$[\ma W]_{\cdot,1} \leftarrow \ve 1 /[\ma R]_{1,1}$\;
	\For{$\ell=2,3, \ldots, |\set I|$}{
		Pick smallest $k$ such that $\exists j$ where  $\set I[k] + \ve e_j = \set I[\ell]$\;
		$\ve v_\ell \leftarrow \diag([ [\ve z_1]_j, \ldots, [\ve z_M]_j]) \ve q_j$\;
		$\ve v_\ell \leftarrow \ve v_\ell - [\ma W]_{\cdot,1:\ell-1} [\ma R]_{1:\ell-1,\ell}$\;
		$[\ma W]_{\cdot, \ell} \leftarrow \ve v_\ell/[\ma R]_{\ell,\ell}$\;
	}
\end{algorithm2e}
\vspace*{-1.5em}
\end{minipage}
\caption{Evaluating a Vandermonde with Arnoldi Basis at New Points}
\label{alg:arnoldi_eval}
\end{algorithm}

% Refinement to local optimality
\section{Refinement to Local Optimality\label{sec:refinement}}
In some situations we desire locally optimal rational approximations;
namely, $\ve a^\star$ and $\ve b^\star$ satisfying the first order necessary conditions of
\begin{equation}\label{eq:refine_ls}
	\min_{\ve a, \ve b} \| \ve f(\ve a, \ve b)\|_2,
	\quad \text{where} \quad
	\ve f(\ve a, \ve b) \coloneqq \ve y - \diag(\ma Q\ve b)^{-1}\ma P \ve a.
\end{equation}
Although the best iterate of the Stabilized SK will not satisfy the local optimality conditions,
it frequently provides a good initialization for a nonlinear least squares solver.
There is only one difficulty in applying standard nonlinear least squares algorithms 
to the rational approximation problem~\cref{eq:refine_ls}:
the Jacobian of $\ve f$ is structurally rank-deficient.
For any scalar $\alpha$, $\ve f(\ve a, \ve b) = \ve f(\alpha \ve a, \alpha \ve b)$.
Hence there is one additional degree of freedom if the coefficients $\ve a$ and $\ve b$ are real;
two if these coefficients are complex.
In our implementation we remove this rank deficiency by fixing 
the value of the largest entry in $\ve b$.

An additional concern is that the Jacobian of $\ve f$,
\begin{align}\label{eq:jacobian}
	\ma F(\ve a, \ve b) &=
	\begin{bmatrix} 
		\diag(\ma Q\ve b)^{-1}\ma P &
		-\diag[\diag(\ma Q\ve b)^{-2}\ma P \ve a]\ma Q
	\end{bmatrix},
\end{align}
can be ill-conditioned due to the presence of $\diag(\ma Q\ve b)^{-1}$
much like the SK iteration.
We can partially rectify this by using 
the basis generated by the $\ell$th step of the Stabilized SK iteration.

\section{Numerical Experiments\label{sec:examples}}
Here we compare the Stabilized SK iteration 
to other univariate and multivariate rational approximation algorithms
on a series of test problems from recent literature.

\subsection{Univariate Problems\label{sec:examples:univariate}}
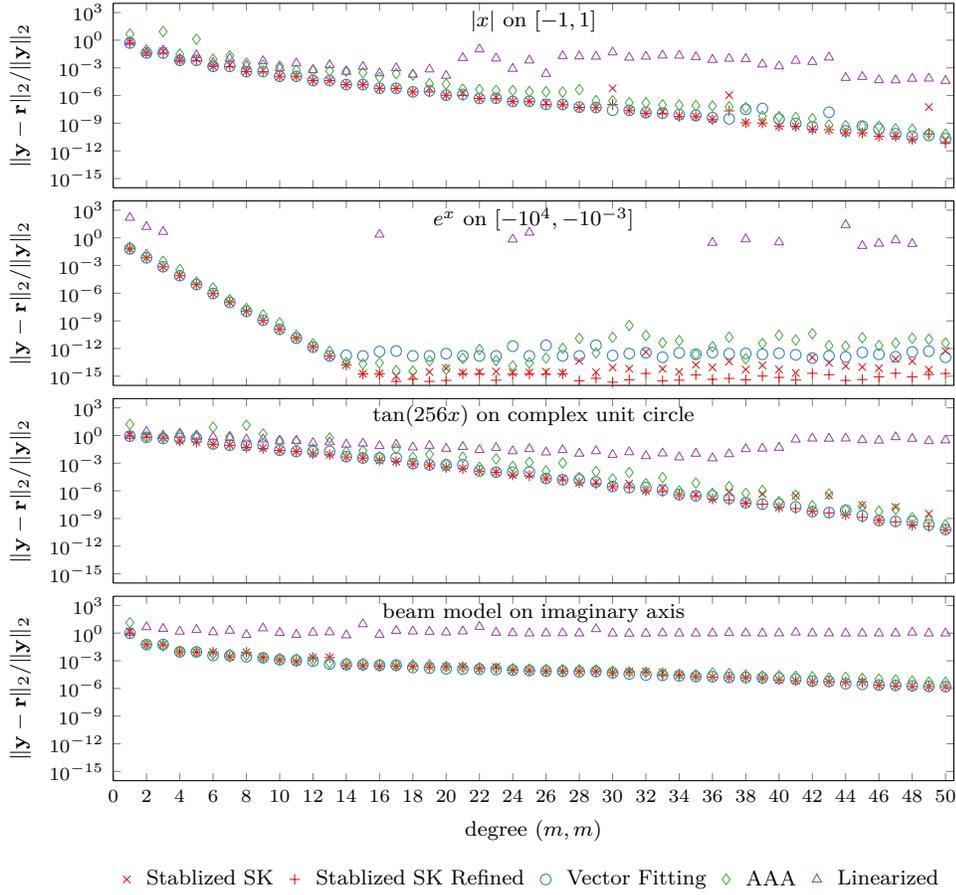
\begin{figure}

\begin{tikzpicture}
\begin{groupplot}[
	group style = {group size = 1 by 4, vertical sep = 5pt, horizontal sep = 7pt},
	title style={yshift=-7pt},
	ymode = log,
	xmin = 0,
	xmax = 50.5,
	ymin = 1e-16,
	ymax = 1e4,
	width = 0.98\linewidth,
	height = 0.31\linewidth,
	xtick = {0,2,4,..., 50},
	ytickten = {-15,-12,...,3},
	every axis plot post/.append style={
		mark size = 2pt,
		},
	xlabel = {degree $(m,m)$},
	ylabel = {$\|\ve y - \ve r\|_2/\|\ve y\|_2$},
	title style={at={(0.5,1)},anchor=north,yshift=2pt},
	]
	\nextgroupplot[
		title = {$|x|$ on $[-1,1]$},
		xticklabels = {,,},
		xlabel = {},
		legend style={
			at={($(0,0)+(1cm,1cm)$)},
			legend columns=5,
			fill=none,
			draw=none,
			anchor=center,
			align=center,
			column sep=4pt
		},
		% https://tex.stackexchange.com/a/272986
		legend image post style={scale=1},
        legend to name=fred
		]

		\coordinate (c1) at (rel axis cs:0,1);

		\pgfplotstableread{data/fig_scalar_abs.dat}\data;
		\addplot[red, only marks, mark=x]
			table [x=m, y = err_isk] {\data};
		\addlegendentry{Stablized SK};    
		\addplot[red, only marks, mark=+]
			table [x=m, y = err_iskr] {\data};
		\addlegendentry{Stablized SK Refined};    
		
		\addplot[blue, only marks, mark=o]
			table [x=m, y = err_vf] {\data};
		\addlegendentry{Vector Fitting};    
		
		\addplot[green, only marks, mark=diamond]
			table [x=m, y = err_aaa] {\data};
		\addlegendentry{AAA};

		\addplot[purple, only marks, mark=triangle]
			table [x=m, y = err_lra] {\data};
		\addlegendentry{Linearized};

	\nextgroupplot[
		title = {$e^x$ on $[-10^4,-10^{-3}]$},
		xticklabels = {,,},
		xlabel = {},
		]

		\coordinate (c2) at (rel axis cs:1,1);
			
		\pgfplotstableread{data/fig_scalar_exp.dat}\data;
		\addplot[red, only marks, mark=x]
			table [x=m, y = err_isk] {\data};
		\addplot[red, only marks, mark=+]
			table [x=m, y = err_iskr] {\data};
		\addplot[blue, only marks, mark=o]
			table [x=m, y = err_vf] {\data};
		
		\addplot[green, only marks, mark=diamond]
			table [x=m, y = err_aaa] {\data};
		
		\addplot[purple, only marks, mark=triangle]
			table [x=m, y = err_lra] {\data};

	\nextgroupplot[
		title = {$\tan(256x)$ on complex unit circle},
		title style={align=center},
		xlabel = {},
		xticklabels = {,,},
		%ytickten = {-4,-3,...,4},
		]
		
		\pgfplotstableread{data/fig_scalar_tan256.dat}\data;
		\addplot[red, only marks, mark=x]
			table [x=m, y = err_isk] {\data};
		\addplot[red, only marks, mark=+]
			table [x=m, y = err_iskr] {\data};
		\addplot[blue, only marks, mark=o]
			table [x=m, y = err_vf] {\data};
		\addplot[green, only marks, mark=diamond]
			table [x=m, y = err_aaa] {\data};
		\addplot[purple, only marks, mark=triangle]
			table [x=m, y = err_lra] {\data};

	\nextgroupplot[
		title = {beam model on imaginary axis},
		%yticklabels = {,,},
		%ylabel = {},
		%xticklabels = {,,},
		%ytickten = {-4,-3,...,4},
		]
		\pgfplotstableread{data/fig_scalar_beam.dat}\data;
		\addplot[red, only marks, mark=x]
			table [x=m, y = err_isk] {\data};
		\addplot[red, only marks, mark=+]
			table [x=m, y = err_iskr] {\data};
		\addplot[blue, only marks, mark=o]
			table [x=m, y = err_vf] {\data};
		\addplot[green, only marks, mark=diamond]
			table [x=m, y = err_aaa] {\data};
		\addplot[purple, only marks, mark=triangle]
			table [x=m, y = err_lra] {\data};

\end{groupplot}
	\coordinate (c3) at ($(c1)!.5!(c2)$);
    \node[below] at (c3 |- current bounding box.south)
      {\pgfplotslegendfromname{fred}};	
\end{tikzpicture}
\caption{The performance of
different rational approximation algorithms 
on the four univariate rational approximation examples described in \cref{sec:examples:univariate}.
Many of the linearized results for the exponential example are outside the plotting range. 
}
\label{fig:scalar}

\end{figure}

Here we consider four univariate rational approximation
test problems from Nakatsukasa, S\`ete, and Trefethen~\cite{NST18}
and compare the performance of 
of AAA~\cite{NST18}, vector fitting~\cite{GS99},
linearized rational approximation~\cite{AKL+19x},
and our Stabilized SK iteration
both with and without refinement to local least-squares optimality.
These results are shown in \cref{fig:scalar}.
In each case we see that the rational approximation generated by the Stabilized SK iteration
yields one of the best approximations with a similar residual norm as to Vector Fitting;
both AAA and the linearized approach yield worse approximations.
We also observe that refining the Stabilized SK approximation to a local optimizer
does not often substantially improve the result.
Details on these four examples follows.

\begin{example}[Absolute Value]
	\label{ex:abs}
	Approximating
	$f(x) = |x|$ using 200,000 equispaced points on the interval $[-1,1]$~\cite[Subsec.~6.7]{NST18}.
\end{example}

This example challenges each algorithm with a large quantity of data.
As the absolute value function is even, 
we only see improvement when the numerator and denominator degrees are even.
Although the Stabilized SK iteration often yields a good rational approximation,
there are some cases where this algorithm fails;
i.e., degrees 30, 37, and 49.
This may be due to extreme ill-conditioning on the first step, as seen in \cref{fig:abs},
preventing the algorithm from making further progress.

\begin{example}[Exponential]
	Approximating
	$f(x) = \exp(x)$ using 2,000 logarithmically spaced points on the interval 
	$[-10^4, -10^{-3}]$~\cite[Subsec.~6.8]{NST18}. 
\end{example}

This example challenges each algorithm to handle points
separated by seven orders of magnitude.
Each algorithm except for the linearized rational approximation performs well in this example.
The linearized approach breaks down because the denominator it identifies
has a near zero-value at every point.

\begin{example}[Tangent]
	Approximating 
	$f(x) = \tan(256x)$ using 1000 equispaced points on the unit circle~\cite[Subsec.~6.3]{NST18}.
\end{example}

This example illustrates behavior using complex-valued points $x_j$ and responses $y_j$.
AAA exhibits some oscillations in the accuracy of its approximation
whereas the Stabilized SK iteration (especially after refinement) converges smoothly.

\begin{example}[Beam]
	Approximating the beam model~\cite{CD02} 
	(a rational function of degree $(347,348)$)
	using 1000 points on the imaginary axis with  
	500 logarithmically spaced points between $10^{-2}i$ and $10^{2}i$
	and their complex conjugates~\cite[Subsec.~6.9]{NST18}.
\end{example}

This example provides a system identification application of rational approximation.
There is no clear best algorithm among Stabilized SK, AAA, and Vector Fitting,
although the overall trend is similar.

\subsection{Spurious Poles\label{sec:examples:spurious}}
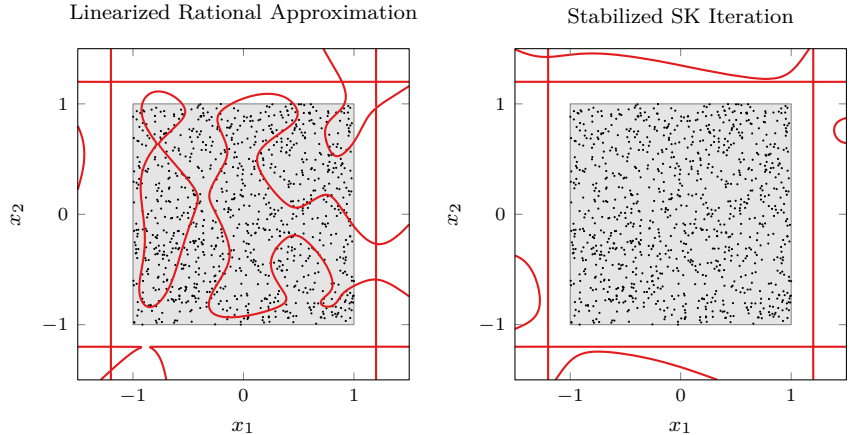
\begin{figure}
\centering

\begin{tikzpicture}
\begin{groupplot}[
		group style = {group size = 2 by 1, horizontal sep = 40pt},
		xmin = -1.5,
		xmax = 1.5,
		ymin = -1.5,
		ymax = 1.5,
		width=0.46\linewidth,
		height=0.46\linewidth,
		axis equal,
		clip mode=individual,
		xlabel = $x_1$,
		ylabel = $x_2$,
	]

	\nextgroupplot[
		title = Linearized Rational Approximation,
	]
	
	\draw[line width=0pt, gray, fill=gray, fill opacity = 0.2] 
		(-1,-1) -- (-1,1) -- (1,1) -- (1,-1) -- cycle;

	\addplot[black, only marks, mark=*, mark size = 0.2pt]
		table [x = x1, y = x2] {data/fig_exp_points.dat};	
	
	\addplot[thick, contour prepared ={draw color=red,labels=false,}, contour prepared format=matlab] 
		table {data/fig_exp_lra.dat};
	
	\nextgroupplot[
		title = Stabilized SK Iteration,
	]
	
	\draw[line width=0pt, gray, fill=gray, fill opacity = 0.2] 
		(-1,-1) -- (-1,1) -- (1,1) -- (1,-1) -- cycle;
	
	\addplot[black, only marks, mark=*, mark size = 0.2pt]
		table [x = x1, y = x2] {data/fig_exp_points.dat};	
	
	\addplot[thick, contour prepared ={draw color=red,labels=false,}, contour prepared format=matlab] 
		table {data/fig_exp_ssk.dat};

\end{groupplot}
\end{tikzpicture}

\caption{The Stabilized SK iteration avoids spurious poles unlike linearized rational approximation.
Here we approximate Example~\ref{ex:exp}
by a rational function of total degree-$(20,20)$
and show zeros of the denominator (lines), 
points used to construct the approximation (dots),
and the region of approximation (shaded square);
cf.~\cite[Fig.~1]{AKL+19x}.
%Note the Stabilized SK iteration avoids spurious zeros in the denominator
%in the region of the approximation 
%while correctly identifying zeros at $x_1= \pm 1.2$ and $x_2 = \pm 1.2$;
}
\label{fig:exp}
\end{figure}

A concern with both AAA~\cite[Sec.~5]{NST18} and the linearized rational approximation~\cite[Subsec.~3.3]{AKL+19x}
is the appearance of Froissart doublets---%
poles of the rational approximation with a near-zero residue
contributing little to the approximation.
Some Froissart doublets are numerical artifacts of the fitting procedure;
others occur even in exact arithmetic as in the absolute value test case (example~\ref{ex:abs})
where every odd degree approximation has one pole with zero residue.
Froissart doublets rarely appear in the rational approximations produced
by the Stabilized SK iteration in our experience.
Fixed points of the SK iteration are often nearby least-squares local minimizers
and these local minimizers are unlikely to have poles with a near-zero residue
as they would contribute little to reducing the residual.
The following example appearing in \cref{fig:exp} illustrates this point.

\begin{example}\label{ex:exp}
	Approximating
	\begin{equation}
		f(x_1, x_2) = \frac{\exp(x_1 x_2)}{(x_1 - 1.2)(x_1 + 1.2)(x_2 - 1.2)(x_2 + 1.2)}
		\quad \text{on} \quad \ve x\in [-1,1]^2
	\end{equation}
	using 1000 randomly distributed points with uniform probability~\cite[Subsec.~3.3]{AKL+19x}.
\end{example}

In this example we expect the rational approximation to be analytic on $[-1,1]^2$
and have poles at $x_1=\pm1.2$ and $x_2 = \pm 1.2$.
In \cref{fig:exp} we see that the linearized approach introduces spurious zeros
inside the approximation domain $[-1,1]^2$ whereas
the Stabilized SK iteration avoids these
for the same degree rational approximation.

\subsection{Parametric Transfer Function Approximation}
One important application of multivariate rational approximation
is parametric model reduction.
In this setting we have a transfer function $H$
that depends on both frequency $z \in \C$
and some parameters $\ve t$, typically real:
\begin{equation}
	H(z, \ve t ) = \ve c^*(z\ma I - \ma A(\ve t))^{-1} \ve b.
\end{equation}
In this context we seek a max-degree rational approximations
as the degree in $z$ is typically much higher than in the parameters $\ve t$.
Here we consider two variants of the Penzl model~\cite[Ex.~3]{Pen06}:
one with one parameter, the other with two.

\begin{example}[One Parameter Penzl Model]\label{ex:penzl1}
	Consider a one parameter variant of the Penzl model~\cite[Subsec.~5.2]{IA14}
	\begin{align}\label{eq:penzl1}
		H(z, t) &= 
			\ve c^\trans \left[
				 z\ma I - 
				\diag(\ma A_1(t), \ma A_2, \ma A_3, \ma A_4)
			\right]^{-1}\ve b \quad \text{with}\\
		\ma A_1(t) &= \begin{bmatrix}
			-1 & t \\ -t & -1
		\end{bmatrix}, \
		\ma A_2 = \begin{bmatrix}
			-1 & 200 \\ -200 & -1
		\end{bmatrix}, \
		\ma A_3 = \begin{bmatrix}
			-1 & 400 \\ -400 & -1
		\end{bmatrix}, 
	\end{align}
	$\ma A_4 = -\diag(1,2,\ldots, 1000)$, 
	and $\ve b= \ve c = [\overbrace{10\cdots10}^6 \overbrace{1\cdots1}^{1000}]$.
	We seek to approximate $H$ where $z\in [0.1, 1000]i$ and $t \in [10,100]$
	using a tensor product grid with 
	100 logarithmically spaced points in $z$ and 30 uniformly spaced points in $t$~\cite[Subsec.~3.2.3]{CG20x}. 
\end{example}

\begin{figure}
\centering
\begin{tikzpicture}
\begin{groupplot}[
		group style = {group size = 1 by 3, vertical sep = 8pt,},
		width=0.9\linewidth,
		height=0.3\linewidth,
		%axis equal,
		clip mode=individual,
		xlabel = $z$,
		ylabel = $t$,
		xmin = 1e-1,
		xmax = 1e3,
		ymin = 10,
		ymax = 100, 
		enlargelimits=false,
		xtickten = {-1,0,1,2,3},
		xmode = log,
		ytick = {10, 50, 100},
		point meta min = -8,
		point meta max = 0,
		colormap name=Blu8,
		title style={at={(0,.5)},anchor=south,rotate =90, yshift=10pt, xshift=-5pt},
		yticklabel pos=right
	]

	\nextgroupplot[title = Linearized,
		xticklabels = {,,},
		xlabel = {},
	]

	\addplot[contour prepared filled= {labels = {false}}] 
		table {data/fig_penzl_contour_lra.dat};

	\nextgroupplot[title = Parametric-AAA,
		xticklabels = {,,},
		xlabel = {},
		]
	
	\addplot[contour prepared filled= {labels = {false}}] 
		table {data/fig_penzl_contour_paaa.dat};

	\addplot[only marks, mark=*, red, mark size = 1.25pt]
		table [x= z, y = s] {data/fig_penzl_paaa_interp.dat};
	
	\nextgroupplot[title = Stabilized SK,
		colorbar horizontal,
		colorbar style = {
				at = {(0, -0.35), anchor=north west},
				width = \pgfkeysvalueof{/pgfplots/parent axis width},
				xlabel = approximation error $|y_j - r_j|$,
				xmin = -8, xmax = 0,
				xtick = {-8, -7, ...,0},
				xticklabels = {$10^{-8}$, $10^{-7}$, $10^{-6}$, 
					$10^{-5}$, $10^{-4}$, $10^{-3}$, $10^{-2}$, $10^{-1}$, $10^0$},
			},
	]
	
	\addplot[contour prepared filled= {labels = {false}}] 
		table {data/fig_penzl_contour_ssk.dat};

\end{groupplot}
\end{tikzpicture}
\caption{The error in rational approximations of 
of total degree $([8,8],[8,8])$
on the one-parameter Penzl model with two variables described in Example~\ref{ex:penzl1}.
Dots denote the p-AAA interpolation points.
The least squares residual norm on the training data 
is $2.203$ for linearized rational approximation, $0.238$ for p-AAA, and $0.0189$ for Stabilized SK. 
}
\label{fig:penzl}
\end{figure}
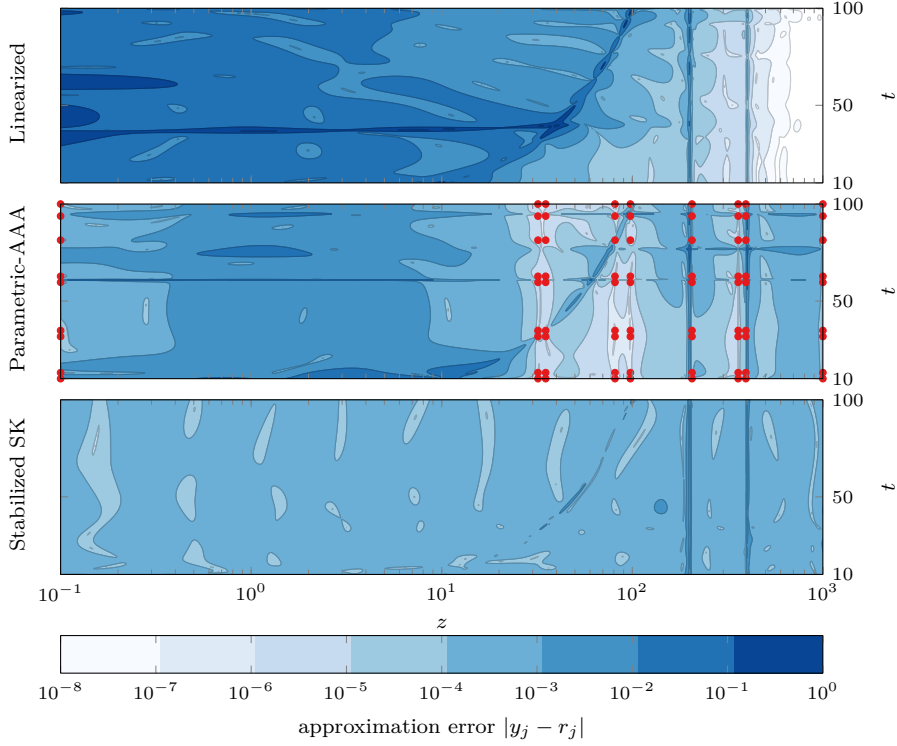

\Cref{fig:penzl} shows the point-wise error of the multivariate rational approximations
produced by linearized rational approximation~\cite{AKL+19x}, Parametric-AAA~\cite{CG20x},
and our Stabilized SK iteration. 
In this case the Stabilized SK iteration produces an approximation
that is accurate throughout the domain
whereas the p-AAA approximation is most accurate near its interpolation points
and the linearized approximation is only accurate for large $z$.
Note the least squares residual norm of Stabilized SK is approximately one tenth that of p-AAA
and approximately one hundredth that of the linearized approach. 

\pagebreak[4]

\begin{example}[Two Parameter Penzl Model]\label{ex:penzl2}
	Consider a two parameter variant of the Penzl model~\cite[Subsec.~3.3.1]{CG20x} 
	\begin{align}
		H(z, t, u) &= 
			\ve c^\trans \left[
				 z\ma I - 
				\diag(\ma A_1(t), \ma A_2(u), \ma A_3(u), \ma A_4)
			\right]^{-1}\ve b \text{ with}\\
		\ma A_2(u) &= \begin{bmatrix}
			-1 & u \\ -u & -1
		\end{bmatrix}, \
		\ma A_3(u) = \begin{bmatrix}
			-1 & 2u \\ -2u & -1
		\end{bmatrix}
	\end{align}
	with the remainder of the variables are as defined in Example~\ref{ex:penzl1}.
	We approximate $H$ for $z \in [1, 2000]i$, $t\in [10,100]$, and $u\in [150,250]$
	using a tensor product grid with $100$ logarithmically spaced frequencies in $z$
	and $10$ equispaced points in $t$ and $u$.
\end{example}

\begin{table}
\caption{A performance comparison of linearized rational approximation, Parametric-AAA, and Stabilized SK 
on the two parameter Penzl model (example~\ref{ex:penzl2})
using the total degrees generated by iterates of p-AAA.}
\label{tab:penzl}
\centering
\begin{filecontents*}{data/fig_penzl3.dat}
d1	d2	d3	paaa_err	lra_err	ssk_err	
6.0	6.0	4.0	0.048494967782378634	0.017457056566102528	0.0010519052884821016	
7.0	6.0	5.0	0.017257131216801506	0.015664943428836356	0.0002464170397436129	
8.0	7.0	5.0	0.0010155676445713579	0.005528321867551583	0.0005086133645922092	
9.0	7.0	6.0	0.0002164115862106732	0.0010274286723870648	9.355833364019876e-05	
10.0	7.0	6.0	2.6579687736374032e-05	0.0006917324300501776	7.215470718325853e-07	
11.0	7.0	6.0	2.853573281060498e-06	0.0006159355456973319	2.029115125375154e-07	
12.0	8.0	7.0	5.008671952829466e-06	0.005010850139860919	1.7921006108176127e-08	
\end{filecontents*}
\pgfplotstabletypeset[
	clear infinite,
	every head row/.style={
		before row={
			\toprule 
			%\multicolumn{6}{c}{degree}  \\
			\multicolumn{3}{c}{numerator} & \multicolumn{3}{c}{denominator} & 
			\multicolumn{6}{c}{$\|\ve y - \ve r\|_2/\|\ve y \|_2$}\\
			\cmidrule(lr){1-3} \cmidrule(lr){4-6} \cmidrule(lr){7-12}
		},
		after row=\midrule,	
	},
	every last row/.style={
		after row=\bottomrule,
	},
	columns/d1/.style = {column name=$z$, fixed},
	columns/d2/.style = {column name=$t$, fixed},
	columns/d3/.style = {column name=$u$, fixed},
	columns/lra_err/.style = {column name = Linearized, sci, precision=4, sci zerofill, dec sep align},
	columns/paaa_err/.style = {column name = Parametric-AAA, sci, precision=4, sci zerofill, dec sep align},
	columns/ssk_err/.style = {column name = Stabilized SK, sci, precision=4, sci zerofill, dec sep align},
	columns = {d1, d2, d3, d1, d2, d3, lra_err, paaa_err, ssk_err},
]{data/fig_penzl3.dat}

\end{table}

\Cref{tab:penzl} shows the history of rational approximations generated by Parametric-AAA
for the two parameter Penzl model
and the corresponding residual norms for both
linearized rational approximation and the  Stabilized SK iteration. 
This example shows the Stabilized SK iteration produces an approximation 
with a smaller residual norm, often by an order of magnitude or more.

\section{Discussion}
% Summary
Here we corrected the numerical instability
in the Sanathanan-Koerner iteration
exposing a practical algorithm for 
univariate and multivariate rational approximations.
Although not optimal, this algorithm yields excellent rational approximations.
Here we have only considered scalar-valued rational approximation problems,
however we could extend the Stabilized Sanathanan-Koerner iteration 
to vector- and matrix-valued
rational approximation following~\cite[subsec.~2.4]{DGB15a}.

%Sometimes starting with the initial weight $\ve 1$
%can yield the first step so ill-conditioned
%that subsequent steps cannot recover.
%A better initialization could fix this,
%which we leave as a topic for future research.

%In some situations is more desirable to rapidly construct a rational approximation
%than to construct a near-optimal approximation; e.g., AAA.
%We could envision increasing the degree of the approximation 
%
%
%There may be better pairs of ordering and update rules;
%for example for total degree polynomials we could choose 
%to multiply by the most recent column rather than the least recent column.

\section*{Acknowledgements}
I would like to thank Caleb Magruder for many helpful discussions in preparing this manuscript
and Yuji Nakatsukasa for the suggestion
to use classical Gram-Schmidt with iterative refinement
in my implementation of multivariate Vandermonde with Arnoldi.

\bibliographystyle{siamplain}
\bibliography{abbrevjournals,master}
\end{document}